# Polyexponentials

Khristo N. Boyadzhiev
Ohio Northern University
Department of Mathematics
Ada, OH 45810
k-boyadzhiev@onu.edu

## 1. Introduction.

The polylogarithmic function [15]

$$\mathrm{Li}_s(x) = \sum_{n=1}^{\infty} \frac{x^n}{n^s} \tag{1.1}$$

and the more general Lerch Transcendent (or Lerch zeta function) [4], [14]

$$\Phi(x, s, \lambda) = \sum_{n=0}^{\infty} \frac{x^n}{(n+\lambda)^s} \tag{1.2}$$

have established themselves as very useful special functions in mathematics. Here and further we assume that $\mathrm{Re}\,\lambda > 0$.

In this article we shall discuss the function

$$e_s(x, \lambda) = \sum_{n=0}^{\infty} \frac{x^n}{n!\,(n+\lambda)^s} \tag{1.3}$$

and survey some of its properties. Our intention is to show that this function can also be useful in a number of situations. It extends the natural exponential function and also the exponential integral (see Section 2). As we shall see, the function $e_s(x, \lambda)$ is relevant also to the theory of the Lerch Transcendent (1.2 ) and the Hurwitz zeta function

$$\zeta(s, \lambda) = \sum_{n=0}^{\infty} \frac{1}{(n+\lambda)^s} \ , \ \mathrm{Re}\,s > 1 \tag{1.4}$$

About a century ago, in 1905, Hardy [6] published a paper on the function

$$F_{a,s}(x) = \sum_{n=0}^{\infty} \frac{x^n}{n!\,(n+a)^s}\,, \qquad\qquad (1.5)$$

where he focused on the variable $x$ and studied the zeroes of $F_{a,s}(x)$ and its asymptotic behavior. We shall use a different notation here in order to emphasize the relation of $e_s(x,\lambda)$ to the natural exponential function. Hardy wrote about $F_{a,s}$: "All the functions which have been considered so far are in many ways analogous to the ordinary exponential function" [6, p. 428]. For this reason we call $e_s(x,\lambda)$ *polyexponential* (or *polyexponentials*, if $s$ is considered a parameter). Our aim is to study some of its properties (complimentary to those in [6]} and point out some applications.

In Section 2 we collect some basic properties of the polyexponential function. Next, in Section 3, we discuss exponential polynomials and the asymptotic expansion of $e_s(x,\lambda)$ in the variable $\lambda$. In Section 4 we evaluate $e_m(x,\lambda)$ for negative integers $m$. Section 5 deals with the Laplace transform of polyexponentials and related integral representations for the Riemann and Hurwitz zeta functions.

An important application is the evaluation of integrals of the form

$$\int_{c-i\infty}^{c+i\infty} x^{-s} R(s)\Gamma(s)\,ds\,, \qquad\qquad (1.6)$$

where $R(s)$ is a rational function. In Section 6 we show how such integrals can be evaluated as combinations of polyexponentials and expressions of the form $e^x p(x) + q(x)$, where $p, q$ are polynomials. In Sections 7 and 8 we evaluate some special series in terms of polyexponentials

Interesting use of $e_s(x,\lambda)$ has been made by Shidlovskii [12], who established algebraic independence. His result is mentioned in the last section.

In what follows we set for brevity

$$e_s(x) = e_s(x,1) = \sum_{n=0}^{\infty} \frac{x^n}{n!\,(n+1)^s}\,. \qquad\qquad (1.7)$$

## 2. Basic properties.

Directly from the definition one has:

$$e_0(x, \lambda) = e^x,$$ 

(2.1)

$$e_1(x) = \frac{e^x - 1}{x},$$

(2.2)

For $x \geq 0$:

$$x^\lambda e_1(-x, \lambda) = \gamma(\lambda, x),$$

(2.3)

where

$$\gamma(\lambda, x) = \int_0^x t^{\lambda-1} e^{-t} dt$$

(2.4)

is the (lower) incomplete Gamma function [1], [10]. This follows from the equation

$$\frac{d}{dx}[x^\lambda e_1(-x, \lambda)] = \sum_{n=0}^\infty \frac{(-1)^n x^{n+\lambda-1}}{n!} = x^{\lambda-1} e^{-x}.$$

(2.5)

When $\mathrm{Re}\,\lambda > 0$ the identity (2.3) extends to all $x \in \mathbb{C} \setminus (-\infty, 0]$ by analytic continuation. When $\lambda$ is a positive integer, (2.3) is true for all $x$.

Next,

$$e_2(x) = \frac{1}{x} \int_0^x \frac{e^t - 1}{t} dt,$$

(2.6)

as

$$(\frac{d}{dx} x) e_2(x) = e_1(x).$$

(2.7)

Equation (2.6) can be written also in the form

$$e_2(x) = \frac{-1}{x} \mathrm{Ein}(-x),$$

(2.8)

or

$$\mathrm{Ein}(x) = x e_2(-x),$$

(2.9)

where $\mathrm{Ein}(z)$ is the special entire function (complementary exponential integral, [1], [10])

$$\text{Ein}(z) = \int\limits_0^z \frac{1 - e^{-t}}{t} dt = \sum_{k=1}^{\infty} \frac{(-1)^{k-1} z^k}{k!\, k} \, . \tag{2.10}$$

Equation (2.7) extends to the recursive relation:

$$(\frac{d}{dx} x)\, e_{p+1}(x) = e_p(x), \tag{2.11}$$

or

$$e_{p+1}(x) = \frac{1}{x} \int\limits_0^x e_p(t)\, dt\, . \tag{2.12}$$

More generally,

$$(\frac{d}{dx} x)\, [x^{\lambda-1}\, e_{p+1}(x, \lambda)] = x^{\lambda-1} e_p(x, \lambda), \tag{2.13}$$

and solving for $e_{p+1}$:

$$e_{p+1}(x, \lambda) = x^{-\lambda} \int\limits_0^x t^{\lambda-1} e_p(t, \lambda)\, dt\, . \tag{2.14}$$

For $p \in \mathbb{N}$ , we obtain from (2.13)

$$(\frac{d}{dx} x)^p\, [x^{\lambda-1} e_p(x, \lambda)] = x^{\lambda-1} e_0(x, \lambda) = x^{\lambda-1} e^x, \tag{2.15}$$

in particular

$$(\frac{d}{dx} x)^p\, e_p(x) = e_0(x) = e^x\, .$$

The next property involves a special infinite series. Let $|z| < |\lambda|$ . Then

$$\sum_{p=0}^{\infty} e_p(x, \lambda) z^p = e^x + z e_1(x, \lambda - z)\, . \tag{2.16}$$

The proof is straightforward - we reverse the order of summation in

$$\sum_{p=0}^{\infty} \sum_{n=0}^{\infty} \frac{x^n z^p}{n!\,(n+\lambda)^p} \tag{2.17}$$

and sum the geometric series with ratio $z/(n+\lambda)$.

**Taylor expansion according to the variable $\lambda$**.

The series (2.16) is in fact a Taylor expansion, since we have

$$(\frac{d}{dz})^m e_1(x,\lambda-z) = m!\,e_{1+m}(x,\lambda-z), \tag{2.18}$$

and (2.16) can be written in the form

$$e_1(x,\lambda-z) = \sum_{m=0}^{\infty} e_{1+m}(x,\lambda)\,z^m. \tag{2.19}$$

**Proposition 2.1**. The function $e_s(x,\lambda-z)$ has the following Taylor series on the powers of

$z$

$$e_s(x,\lambda-z) = \sum_{m=0}^{\infty} \frac{\Gamma(s+m)}{\Gamma(s)\,m!} e_{s+m}(x,\lambda)\,z^m, \tag{2.20}$$

similar to the Taylor expansion of the Lerch zeta function - see [7], [13].

*Proof*. Equation (2.18) extends to

$$(\frac{d}{dz})^m e_s(x,\lambda-z) = \frac{\Gamma(s+m)}{\Gamma(s)\,m!} e_{s+m}(x,\lambda-z), \tag{2.21}$$

and the representation follows.

**Hankel integral representation**.

Consider the (Hankel) contour $L$ consisting of three parts: $L = L_- \cup L_+ \cup L_\epsilon$, with $L_-$ the

"lower side" (i.e. $arg(z) = -\pi$) of the ray $(-\infty, -\epsilon)$, $\epsilon > 0$, traced left to right, and $L_+$ the

"upper side" ($arg(z) = \pi$) of this ray traced right to left. Finally, $L_\epsilon = \{z = \epsilon e^{\theta i} : -\pi \le \theta \le \pi\}$ is a small circle traced counterclockwise and connecting the two sides of the ray. This contour is used, for example, in [10, p. 37-38] and its symmetric version extending to $+\infty$ can be found in [4,

p. 25].

**Proposition 2.2.** For all $\operatorname{Re} \lambda > 0, x \in \mathbb{C}$, the polyexponential function has the integral representation

$$e_s(x,\lambda) = \frac{\Gamma(1-s)}{2\pi i} \int_L z^{s-1} e^{\lambda z} e^{xe^z} \, dz, \tag{2.22}$$

when $s \notin \mathbb{N}$, and

$$e_m(x,\lambda) = \frac{(-1)^m}{2\pi i (m-1)!} \int_L z^{m-1} e^{\lambda z} e^{xe^z} \operatorname{Log} z \, dz, \tag{2.23}$$

when $m \in \mathbb{N}$.

*Proof.* Starting from the well-known representation [10, (1.12), p.38]

$$\frac{1}{\Gamma(1-s)} = \frac{1}{2\pi i} \int_L z^{s-1} e^z \, dz, \tag{2.24}$$

$(s \notin \mathbb{N})$ we derive

$$\frac{1}{(n+\lambda)^s} = \frac{\Gamma(1-s)}{2\pi i} \int_L z^{s-1} e^{\lambda z} e^{nz} \, dz, \tag{2.23}$$

and (2.22) follows from here - multiply both sides by $x^n/n!$ and sum for $n = 0, 1, \ldots$ . Next, in order to prove (2.23) we write first

$$\Gamma(1-s) = \frac{\pi}{\Gamma(s)\sin \pi s},$$

and then use the rule of L'Hospital to evaluate the limit

$$e_m(x,\lambda) = \frac{1}{2\pi i} \lim_{s \to m} \frac{\pi}{\Gamma(s)\sin \pi s} \int_L z^{s-1} e^{\lambda z} e^{xe^z} \, dz$$

$$= \frac{1}{2\pi i (m-1)!} \lim_{s \to m} \frac{\pi}{\sin \pi s} \int_L z^{s-1} e^{\lambda z} e^{xe^z} \, dz,$$

thus obtaining (2.23).

**Asymptotic for large $x$.** For completeness we include here also the asymptotic for real $x \to \pm\infty$ as obtained by Hardy [6]. For sufficiently large $x$, there exists a function $\epsilon_x$ tending uniformly to zero like $1/x$ when $x \to +\infty$ such that:

$$e_s(x, \lambda) = \frac{e^x}{x}(1 + \epsilon_x) \tag{2.26}$$

$$e_s(-x, \lambda) = \frac{\Gamma(\lambda)}{\Gamma(s)}\frac{(\log x)^{s-1}}{x^\lambda}(1 + \epsilon_x). \tag{2.27}$$

**Symbolic representation**. The function $e_s(x, \lambda)$ can be associated with the Mellin fractional derivative

$$(x\frac{d}{dx} + \lambda)^\alpha, \tag{2.28}$$

where $\alpha \in \mathbb{C}$. The action of this derivative on a given analytical function $f(x) = \sum_{n=0}^\infty a_n x^n$ is described by the equation

$$(x\frac{d}{dx} + \lambda)^\alpha f(x) = \sum_{n=0}^\infty a_n(n + \lambda)^\alpha x^n. \tag{2.29}$$

Thus:

$$e_s(x, \lambda) = (x\frac{d}{dx} + \lambda)^{-s} e^x. \tag{2.30}$$

The analytical theory of the derivative (2.28) is developed in [3].

## 3 Asymptotic expansion in $\lambda$ and exponential polynomials.

The asymptotic expansion of $e_s(x, \lambda)$ in the variable $\lambda$ was obtained in [2]. Namely,

$$e_s(x, \lambda) = \sum_{n=0}^\infty \frac{x^n}{n!(n + \lambda)^s} \sim e^x \sum_{n=0}^\infty \binom{-s}{n}\frac{1}{\lambda^{n+s}}\phi_n(x), \tag{3.1}$$

where

$$\phi_n(x) = e^{-x}(x\frac{d}{dx})^n e^x = e^{-x}\sum_{k=0}^\infty \frac{k^n}{k!}x^k, \quad n = 0, 1, \dots \tag{3.2}$$

are the exponential polynomials [2]. We list here some of their properties:

$$\phi_{n+1} = x \left( \phi_n' + \phi_n \right), \tag{3.3}$$

$$\phi_{n+1} = x \sum_{k=0}^{n} \binom{n}{k} \phi_k, \tag{3.4}$$

$$\phi_n' = \sum_{k=0}^{n-1} \binom{n}{k} \phi_k, \tag{3.5}$$

$$\phi_0(x) = 1, \tag{3.6}$$

$$\phi_1(x) = x, \tag{3.7}$$

$$\phi_2(x) = x + x^2, \tag{3.8}$$

$$\phi_3(x) = x + 3x^2 + x^3, \text{ etc.} \tag{3.9}$$

In general

$$\phi_n(x) = \sum_{k=0}^{n} \left\{ {n \atop k} \right\} x^k, \tag{3.10}$$

where $\left\{ {n \atop k} \right\}$ are the Stirling numbers of the second kind (the number of partitions of a set of

$n$ elements into $k$ disjoint nonempty subsets [5]).

The exponential generating function of these polynomials is

$$e^{x(e^t - 1)} = \sum_{n=0}^{\infty} \phi_n(x) \frac{t^n}{n!}. \tag{3.11}$$

We shall need also the following lemma.

**Lemma 3.1**. For $p = 0, 1, 2, \ldots$:

$$\int_0^\infty e^{-2t} \phi_p(-t) \, dt = \frac{1 - 2^{p+1}}{p+1} B_{p+1}, \tag{3.12}$$

where $B_p$ are the Bernoulli numbers. This follows from (3.10) and the identity [5, p.317, Problem

6.76]:

$$\sum_{k=0}^{p} \left\{{p \atop k}\right\} k! \frac{(-1)^k}{2^{k+1}} = (1 - 2^{p+1}) \frac{B_{p+1}}{p+1} \tag{3.13}$$

Another interesting property involving $\phi_k$ and the Bernoulli numbers is (8.6) in section 8 .

### 4 Values at negative integer $s$ in terms of exponential polynomials.

It is known that for $n = 1, 2, \ldots,$

$$\zeta(-n, \lambda) = \frac{-B_{n+1}(\lambda)}{n+1}, \tag{4.1}$$

where $B_k(\lambda)$ are the Bernoulli polynomials [1], [4]. Similarly, for $p \in \mathbb{N}$, $e_{-p}(x, \lambda) e^{-x}$ are polynomials in $\lambda, x$. This was briefly mentioned by Hardy [6] on pp. 426-427, where he found a recurrent relation for these polynomials. We shall list them here in explicit form.

Let $s = -p$, where $p = 0, 1, \ldots$. Then for every $x$ and every $\mathrm{Re}\,\lambda > 0$,

$$e_{-p}(x, \lambda) = \sum_{n=0}^{\infty} (n + \lambda)^p \frac{x^n}{n!} \tag{4.2}$$

$$= \sum_{n=0}^{\infty} \left[ \sum_{k=0}^{p} \binom{p}{k} n^k \lambda^{p-k} \right] \frac{x^n}{n!} = \sum_{k=0}^{p} \binom{p}{k} \lambda^{p-k} \sum_{n=0}^{\infty} n^k \frac{x^n}{n!},$$

and according to (3.2):

$$e_{-p}(x, \lambda) = e^x \sum_{k=0}^{p} \binom{p}{k} \lambda^{p-k} \phi_k(x). \tag{4.3}$$

In particular, when $\lambda = 1$ we find in view of (3.4):

$$e_{-p}(x) = e^x \sum_{k=0}^{p} \binom{p}{k} \phi_n(x) = \frac{e^x}{x} \phi_{p+1}(x), \tag{4.4}$$

or

$$\phi_p(x) = x e^{-x} e_{1-p}(x). \tag{4.5}$$

Using this equation we can extend $\phi_p$, as a function of $p$, when $p \neq 0, 1, \ldots$, in terms of the polyexponential $e_{1-p}$.

For future reference we denote the polynomials in (4.3) by:

$$Q_p(x, \lambda) = \sum_{k=0}^{p} \binom{p}{k} \lambda^{p-k} \phi_k(x),$$ (4.6)

that is, on the powers of $x$:

$$Q_0(x, \lambda) = 1,$$ (4.7)

$$Q_1(x, \lambda) = x + \lambda,$$ (4.8)

$$Q_2(x, \lambda) = x^2 + (2\lambda + 1) x + \lambda^2,$$ (4.9)

$$Q_3(x, \lambda) = x^3 + (3\lambda + 3) x^2 + (3\lambda^2 + 3\lambda + 1) x + \lambda^3, \text{ etc.}$$ (4.10)

## 5. Laplace transform and relation to the Riemann zeta function.

We evaluate here the Laplace transform of the polyexponential function in the variable $x$. The Laplace image will be written in the variable $y > 0$,

$$\mathsf{L}[e_s(x, \lambda)](y) = \int_0^{\infty} e_s(x, \lambda) e^{-xy} dx$$ (5.1)

$$= \sum_{n=0}^{\infty} \frac{1}{n!(n+\lambda)^s} \int_0^{\infty} x^n e^{-xy} dx = \sum_{n=0}^{\infty} \frac{1}{(n+\lambda)^s} \frac{1}{y^{n+1}},$$

that is,

$$\mathsf{L}[e_s(x, \lambda)](y) = \frac{1}{y} \Phi\left(\frac{1}{y}, s, \lambda\right),$$ (5.2)

where $\Phi$ is the Lerch Transcendent (1.2). In particular, for $\mathbf{Re}\, s > 1$, $y = 1$:

$$\mathsf{L}[e_s(x, \lambda)](1) = \zeta(s, \lambda).$$ (5.3)

Equations (5.2) and (5.3) can be written also in the form

$$\Phi(x,s,\lambda) = \int\limits_0^\infty e_s(tx,\lambda)e^{-t}dt, \qquad (5.4)$$

and correspondingly,

$$\zeta(s,\lambda) = \int\limits_0^\infty e_s(t,\lambda)e^{-t}dt. \qquad (5.5)$$

In particular, when $\lambda = 1$,

$$\zeta(s) = \int\limits_0^\infty e_s(t)e^{-t}dt \qquad (5.6)$$

Note that (5.4) is also true for $\operatorname{Re} s \le 1$ when $|x| < 1$. The representations (5.5) and (5.6), however, do not extend to $\operatorname{Re} s \le 1$. For example, when $s = -p, p = 0, 1, \ldots$, we have by (4.4)

$$e_{-p}(t) = \frac{e^t}{t}\phi_{p+1}(t),$$

and the integral in (5.6) is divergent. We shall modify now equations (5.5) and (5.6) in order to obtain representations valid for all complex $s$. Replacing $e_s(t,\lambda)$ by $e_s(-t,\lambda)$ we compute:

$$\int\limits_0^\infty e_s(-t)e^{-t}dt = \sum_{n=0}^\infty \frac{(-1)^n}{n!(n+\lambda)^s}\int\limits_0^\infty t^n e^{-t}dt = \sum_{n=0}^\infty \frac{(-1)^n}{(n+\lambda)^s}. \qquad (5.7)$$

The series on the right side converges for all $\operatorname{Re} s > 0$. This function is a special case of the Lerch zeta function, $\Phi(-1,s,\lambda)$. We shall use the notation:

$$\eta(s,\lambda) = \sum_{n=0}^\infty \frac{(-1)^n}{(n+\lambda)^s}. \qquad (5.8)$$

**Proposition 5.1**. For all complex $s$:

$$\eta(s,\lambda) = \int\limits_0^\infty e_s(-t,\lambda)e^{-t}dt. \qquad (5.9)$$

*Proof*. The function $\eta(s,\lambda)$ extends to an entire function in the variable $s$ (see[14]). The integral in (5.9) also extends as a function defined for all $s$. This can be seen by integrating (2.18) and changing the order of integration. Thus we obtain for $\mathrm{Re}\,s > 0$, $s \neq 1, 2, 3, \ldots,$

$$\eta(s, \lambda) = \int_0^\infty e_s(-t, \lambda) e^{-t} dt = \frac{\Gamma(1-s)}{2\pi i} \int_L \frac{z^{s-1} e^{\lambda z}}{1 + e^z} dz \ , \tag{5.10}$$

and since the Hankel integral exists for all $s$, equality for all $s$ in (5.10) follows by analytic continuation.

In fact, the representation (5.9) independently shows that $\eta(s, \lambda)$ extends to an entire function of $s$.

When $\lambda = 1$,

$$\eta(s, 1) = \eta(s) = \sum_{k=1}^\infty \frac{(-1)^{k-1}}{k^s}, \tag{5.11}$$

is the Dirichlet eta function with $\eta(s) = \zeta(s)(1 - 2^{1-s})$. Therefore, we have:

**Corollary 5.2.** For all complex $s$,

$$\zeta(s)(1 - 2^{1-s}) = \int_0^\infty e_s(-t) e^{-t} dt \ . \tag{5.12}$$

For example, setting $s = -p$, $p \in \mathbb{N}$, in (5.12) we find

$$\zeta(-p)(1 - 2^{p+1}) = \int_0^\infty e^{-2t}(\phi_p(-t) + \phi_p'(-t)) dt = \int_0^\infty e^{-2t} \phi_p(-t) dt \ ,$$

and in view of (3.12) we arrive at the classical Euler result:

$$\zeta(-p) = \frac{-B_{p+1}}{p+1}. \tag{5.13}$$

It is known also that for all $p \in \mathbb{N}$:

$$\eta(-p, \lambda) = \frac{1}{2} E_p(\lambda),$$

where $E_p(\lambda)$ are the Euler polynomials [14] . Combining (3.12), (4.3) and (5.9) we obtain

$$E_p(\lambda) = 2 \sum_{k=0}^{p} \binom{p}{k} \lambda^{p-k} \frac{1 - 2^{k+1}}{k+1} B_{k+1}, \qquad (5.14)$$

or, in a different form:

$$\eta(-p, \lambda) = -\sum_{k=0}^{p} \binom{p}{k} \lambda^{p-k} \eta(-k). \qquad (5.15)$$

### 6. Mellin transform and evaluation of Mellin integrals.

The Mellin transform of one function $g(x)$ is given by

$$G(s) = \int_0^\infty x^{s-1} g(x) dx, \qquad (6.1)$$

(see [9]) with inversion

$$g(x) = \frac{1}{2\pi i} \int_{(c)} x^{-s} G(s) ds, \qquad (6.2)$$

where $c$ is some appropriate real number and the notation $(c)$ stands for integration upward on the vertical line with abscissa $c$, as in (1.6). Popular examples are

$$\Gamma(s) = \int_0^\infty x^{s-1} e^{-x} dx, \qquad (6.3)$$

and

$$e^{-x} = \frac{1}{2\pi i} \int_{(c)} x^{-s} \Gamma(s) ds \quad (c > 0). \qquad (6.4)$$

**Proposition 6.1**. The Mellin transform of $e_p(-x, \lambda)$ in the $x$- variable is

$$\frac{\Gamma(s)}{(\lambda - s)^p} = \int_0^\infty x^{s-1} e_p(-x, \lambda) dx \quad (0 < \mathrm{Re}\, s < \mathrm{Re}\, \lambda), \qquad (6.5)$$

i.e.

$$e_p(-x, \lambda) = \frac{1}{2\pi i} \int_{(c)} x^{-s} \frac{\Gamma(s)}{(\lambda - s)^p} ds \quad (0 < c < \mathrm{Re}\, \lambda). \qquad (6.6)$$

*Proof.* (6.6) is verified by closing the contour of integration to the left and evaluating the integral by the residue theorem. We use the fact that $\operatorname*{Res}_{s=-n} \Gamma(s) = (-1)^n/n!$.

One interesting corollary from (6.5) is obtained by taking $s = 1/2$ and $\lambda = 1$. Thus we find the evaluation

$$\int_0^\infty \frac{e_p(-x)}{\sqrt{x}} dx = 2^p \sqrt{\pi} \quad (p \geq 0). \tag{6.7}$$

To understand better the singularity at $s = \lambda$ on the left side in (6.5), one may look at the asymptotic (2.26). At $s = \lambda$ the integral becomes divergent.

Replacing $p$ by $-p$, $p = 0, 1, \ldots$ in (6.5) we arrive at the representation

$$\Gamma(s)(\lambda - s)^p = \int_0^\infty x^{s-1} e_{-p}(-x, \lambda) dx = \int_0^\infty x^{s-1} e^{-x} Q_p(-x, \lambda) dx, \tag{6.8}$$

where $Q_p$ is the polynomial defined in (4.6). For $s = \lambda$ we obtain the interesting identity

$$\int_0^\infty x^{\lambda-1} e_{-p}(-x, \lambda) dx = \int_0^\infty x^{\lambda-1} e^{-x} Q_p(-x, \lambda) dx = 0, \tag{6.9}$$

true for all $\operatorname{Re} \lambda > 0$ and $p = 1, 2, \ldots$. For instance, when $p = 1$ this is the well-known property

$$\Gamma(\lambda + 1) - \lambda \Gamma(\lambda) = 0. \tag{6.10}$$

We present now a method for the evaluation of some Mellin integrals. First, this lemma:

**Lemma 6.2**. Suppose $f(s) = a_0 + a_1 s + \ldots + a_m s^m$ is a polynomial. Then for any $c > 0$,

$$\frac{1}{2\pi} \int_{(c)} x^{-s} f(-s) \Gamma(s) ds = e^{-x} \sum_{k=0}^m a_k \phi_k(-x), \tag{6.11}$$

*Proof*. Applying the operator $(x \frac{d}{dx})^k$ to both sides in (6.4) we find according to (3.2)

$$e^{-x} \phi_k(-x) = \frac{1}{2\pi} \int_{(c)} x^{-s} (-s)^k \Gamma(s) ds, \tag{6.12}$$

for $k = 0, 1, \ldots, m$. The lemma follows immediately.

**Theorem 6.3**. Let $R(s)$ be a rational function whose poles $\{\lambda_k\}$ are in some half plane of the form $\mathrm{Re}\, s > c > 0$ . Then the integral

$$\frac{1}{2\pi} \int\limits_{(c)} x^{-s} R(s) \Gamma(s)\, ds \qquad (6.13)$$

can be evaluated as a finite combination of polyexponentials plus an expression like (6.11), where $f(-s)$ is the polynomial part of $R(s)$.

The proof follows from (6.6) and the representation by partial fractions

$$R(s) = f(-s) + \sum_{k,j} \frac{A_{kj}}{(\lambda_k - s)^{p_{kj}}} , \qquad (6.14)$$

thus

$$\frac{1}{2\pi} \int\limits_{(c)} x^{-s} R(s) \Gamma(s)\, ds = e^{-x} \sum_{k=0}^{m} a_k \phi_k(-x) + \sum_{k,j} A_{kj}\, e_{p_{kj}}(-x, \lambda_k). \qquad (6.15)$$

In the general case when the poles of the rational function $R(s)$ are on both sides of $c$ , a simple adjustment is needed. We move the line of integration to the left until all poles remain on its right side, collecting the residues at the poles. The value of the integral will be a modification of the above result. Expressions of the form

$$\sum_{n=r}^{\infty} \frac{(-1)^n x^n}{n!\,(n+\lambda)^s} = e_s(-x, \lambda) + polynomial, \qquad (6.16)$$

will also appear.

Finally, we note that the function $\Gamma(s)\, e_s(x, \lambda)$ is itself a Mellin transform in the variable $s$ (cf [6]):

$$\Gamma(s)\, e_s(x, \lambda) = \int\limits_{0}^{\infty} t^{s-1} e^{-\lambda t} e^{x e^{-t}} dt, \qquad (6.17)$$

To prove this we use the expansion

$$e^{x e^{-t}} = \sum_{n=0}^{\infty} \frac{x^n e^{-nt}}{n!} \qquad (6.18)$$

and integrate termwise.

## 7. Some special series.

The following two series are listed as entries 5.2.17 (10) and 5.2.17 (12) in [11, p. 716]

$$\sum_{n=1}^{\infty} \frac{x^n}{n!} \left(1 + \frac{1}{2} + \ldots + \frac{1}{n}\right) = e^x (\gamma + \ln|x| - \mathrm{Ei}(-x)), \tag{7.1}$$

$$\sum_{n=1}^{\infty} \frac{x^n}{n!} (1^3 + 2^3 + \ldots + n^3) = \frac{1}{4} x e^x (x^3 + 8x^2 + 14x + 4) \tag{7.2}$$

Consider now the general series

$$h_s(x, \lambda, w) = \sum_{n=1}^{\infty} \frac{x^n}{n!} \left(\frac{1}{\lambda^s} + \frac{w}{(\lambda+1)^s} + \ldots + \frac{w^{n-1}}{(\lambda+n-1)^s}\right), \tag{7.3}$$

for $\mathrm{Re}\,\lambda > 0$, $|w| \le 1$, $x \in \mathbb{C}$. This series includes (7.1) and (7.2) as particular cases. For instance, $h_s(x, 1, 1)$ is the exponential generating function of the generalized harmonic numbers:

$$H_n^s = \frac{1}{1^s} + \frac{1}{2^s} + \ldots + \frac{1}{n^s}. \tag{7.4}$$

We shall evaluate the series (7.3) in terms of polyexponentials. Considering $h_s(x, \lambda, w)$ as a function of $x$, it is easy to see that

$$\frac{d}{dx} h_s - h_s = e_s(xw, \lambda), \tag{7.5}$$

which can be written as

$$\frac{d}{dx} (e^{-x} h_s) = e^{-x} e_s. \tag{7.6}$$

Integrating we find:

$$h_s(x, \lambda, w) = e^x \int_0^x e^{-t} e_s(tw, \lambda) dt. \tag{7.7}$$

When $s = \lambda = 1$, this is the representation

$$h_1(x) = \sum_{n=1}^{\infty} \frac{x^n}{n!} \left(1 + \frac{w}{2} + \ldots + \frac{w^{n-1}}{n}\right) = \frac{e^x}{w} \left[\text{Ein}(x) - \text{Ein}(x(1-w))\right], \qquad (7.8)$$

which becomes (7.1) for $w = 1$ in view of the relation between the functions Ein and Ei, see [1].

Equation (7.7) provides an interesting representation of the Lerch zeta function. Namely,

$$\Phi(w, s, \lambda) = \lim_{x \to \infty} e^{-x} h_s(x, \lambda, w), \qquad (7.9)$$

is true for $\text{Re}(s) > 1$, $|x| \leq 1$, or all $s$ when $|x| < 1$ according to (5.4). In particular, for $\text{Re}(s) > 1$

$$\zeta(s, \lambda) = \lim_{x \to \infty} e^{-x} h_s(x, \lambda, 1), \qquad (7.10)$$

as follows from (5.5). Also, it follows from (5.9) that for all complex $s$,

$$\eta(s, \lambda) = \lim_{x \to \infty} e^{-x} h_s(x, \lambda, -1). \qquad (7.11)$$

These equation express the fact that the series (1.2), (1.5) and (5.8) are Borel convergent.

**Borel summation.** A series

$$\sum_{n=1}^{\infty} a_n \qquad (7.12)$$

is said to Borel converge to $M$, if

$$M = \lim_{x \to +\infty} e^{-x} \sum_{n=1}^{\infty} (a_1 + a_2 + \ldots + a_n) \frac{x^n}{n!}.$$

If the series (7.12) converges to $M$ in ordinary sense, it also Borel converges to $M$ (see [8, pp. 471-473]). Equations (7.9), (7.10) simply express this fact for the series (1.2) and (1.5). Equation (7.11), however, reveals more than just Borel convergence, since it holds for all complex $s$, while (5.8) is true only for $\text{Re}\, s > 0$.

We continue now with the study of (7.3). Setting for brevity $h_s(x, \lambda) = h_s(x, \lambda, 1)$, i.e.

$$h_s(x,\lambda) = \sum_{n=1}^{\infty} \frac{x^n}{n!} \left( \frac{1}{\lambda^s} + \frac{1}{(\lambda+1)^s} + \ldots + \frac{1}{(\lambda+n-1)^s} \right),$$

we point out three representations for this function. First, from (6.5) we find the Mellin integral representation for $\operatorname{Re} s > 1$:

$$\Gamma(s) h_s(x,\lambda) = \int_0^{\infty} \frac{t^{s-1} e^{-\lambda t}}{1 - e^{-t}} [e^x - e^{xe^{-t}}] dt. \tag{7.13}$$

Another interesting representation is obtained by changing the order of summation:

$$h_s(x,\lambda) = \sum_{n=0}^{\infty} \frac{1}{(n+\lambda)^s} \left\{ \sum_{k=n+1}^{\infty} \frac{x^k}{k!} \right\}. \tag{7.14}$$

Finally, combining (3.1) with (7.7) we find the asymptotic expansion of $h_s(x,\lambda)$ in the variable $\lambda > 0$,

$$h_s(x,\lambda) \sim e^x \sum_{n=0}^{\infty} \left\{ \binom{-s}{n} \int_0^x \phi_n(t)\, dt \right\} \frac{1}{\lambda^{n+s}} \ .$$

## 8. Evaluating $h_{-p}(x, 1, \pm 1)$, $p \in \mathbb{N}$.

For $p = 0, 1, \ldots,$ we have explicitly

$$h_{-p}(x,1,w) = \sum_{n=1}^{\infty} \frac{x^n}{n!} [1^p + w\, 2^p + \ldots + w^{n-1} n^p]. \tag{8.1}$$

With the notation $h_s(s) = h_s(x,1,1)$ for $w = 1$,

$$h_{-p}(x) = \sum_{n=1}^{\infty} \frac{x^n}{n!} [1^p + 2^p + \ldots + n^p]. \tag{8.2}$$

Using (4.3) in combination with (7.7) we obtain the evaluation:

$$h_{-p}(x) = e^x \int_0^x \phi_{p+1}(t) \frac{dt}{t} ,. \tag{8.3}$$

When $p = 3$, this is the series (7.2). In view of (3.10) we can write (8.3) in the form

$$h_{-p}(x) = e^x \sum_{k=1}^{p+1} \{^{p+1}_k\} \frac{x^k}{k},$$ (8.4)

and also in the form

$$h_{-p}(x) = e^x \phi_p(x) + e^x \int_0^x \phi_p(t)\,dt,$$ (8.5)

by means of (3.3).

**Proposition 8.1.** The following identity is true

$$\int_0^x \phi_p(t)\,dt = \frac{1}{p+1} \sum_{k=1}^{p+1} \binom{p+1}{k} B_{p+1-k} \phi_k(x),$$ (8.6)

where $B_k$ are the Bernoulli numbers.

*Proof.* We shall evaluate now $h_{-p}(x)$ in terms of Bernoulli numbers. Starting from the

Faulhaber formula [5, p.283]

$$1^p + 2^p + \ldots + (n-1)^p = \frac{1}{p+1} \sum_{k=0}^{p+1} \binom{p+1}{k} B_{p+1-k} n^k,$$ (8.7)

we compute according to (3.2)

$$\sum_{n=1}^{\infty} \frac{x^n}{n!} [1^p + 2^p + \ldots + (n-1)^p + n^p]$$

$$= \sum_{n=1}^{\infty} n^p \frac{x^n}{n!} + \frac{1}{p+1} \sum_{k=1}^{p+1} \binom{p+1}{k} B_{p+1-k} \sum_{n=1}^{\infty} n^k \frac{x^n}{n!}$$

$$= e^x \phi_p(x) + \frac{e^x}{p+1} \sum_{k=1}^{p+1} \binom{p+1}{k} B_{p+1-k} \phi_k(x).$$ (8.8)

Comparing (8.5) to (8.8) we arrive at (8.6). The proof is completed.

Next we evaluate the series

$$h_{-p}(x, 1, -1) = \sum_{n=1}^{\infty} \frac{x^n}{n!} [1^p - 2^p + \ldots + (-1)^{n-1} n^p], \qquad (8.9)$$

for $p = 0, 1, \ldots$. By (7.7)

$$h_{-p}(x, 1, -1) = e^x \int_0^x e^{-t} e_{-p}(-t) dt \qquad (8.10)$$

$$= e^x \int_0^x e^{-2t} \frac{\phi_{p+1}(-t)}{-t} dt = e^x \int_0^x e^{-2t} [\phi_p'(-t) + \phi_p(-t)] dt,$$

and a simple computation brings this to the compact form

$$\sum_{n=1}^{\infty} \frac{x^n}{n!} [1^p - 2^p + \ldots + (-1)^{n-1} n^p] = -\phi_p(-x) e^{-x} - e^x \int_0^x e^{-2t} \phi_p(-t) dt. \qquad (8.11)$$

### 9. Shidlovskii's theorem

In a different setting, as shown by Shidlovskii [12], the polyexponentials are useful in the theory of transcendental numbers. A classical theorem of Lindemann says that if $\xi \neq 0$ is an algebraic number, then $e^\xi$ is transcendental. This theorem can be extended to finite sets of numbers. For every $n \in \mathbb{N}$, the numbers $a_1, a_2, \ldots, a_n$ are called algebraically independent, if they are not a solution of any algebraic polynomial of $n$ variables. Shidlovskii proved the following theorem [12, 7.3.3]:

**Theorem.** Let $\lambda$ be a non-integer rational number and $\xi \neq 0$ be an algebraic number. Then for every $m \in \mathbb{N}$ the numbers $e^\xi$, $e_p(\xi, \lambda)$, $p = 1, 2, \ldots, m$ are algebraically independent.

**References**

1.   **M. Abramowitz and I.A. Stegun**, *Handbook of mathematical functions*. Dover, New York, 1970

2.   **Khristo N. Boyadzhiev**, A Series Transformation Formula and Related Polynomials, *Int. J. Math. Math. Sci*., 2005:23 (2005) 3849-3866.

3.   **Paul L. Butzer, Anatoly A. Kilbas, Juan J. Trujilo**, Fractional Calculus in the Mellin setting and Hadamard-type fractional integrals, *J. Math. Anal. Appl*., 269, No.1 (2002), 1-27; ibid v. 269, No.2, (2002), 387-400; ibid v. 270, No. 1 (2002), 1-15.

4.   **A. Erdélyi** (editor), *Higher Transcendental Functions*, vol.1, New York: McGrow-Hill,1955.

5.   **L. Graham, Donald E. Knuth, Oren Patashnik**, *Concrete Mathematics*, Addison-Wesley Publ. Co., New York, 1994.

6.   **G. H. Hardy,** On the zeros of certain classes of integral Taylor series, Part II, *Proc. London Math. Soc.* (2), 2 (1905), 401-431. (*Collected Papers*, vol. IV, Clarendon Press, Oxford, 1969).

7.   **Dieter Klusch**, On the Taylor expansion of the Lerch zeta-function, *J. Math. Anal. Appl*.,170 (1992), 513-523.

8.   **Konrad Knopp**, *Theory and Application of Infinite Series*, Dover, New York, 1990.

9.   **F. Oberhettinger**, *Tables of Mellin Transforms*, Springer, 1974.

10.  **Frank W. Olver**, *Asymptotics and Special Functions*, Academic Press, Boston, 1974.

11.  **A. P. Prudnikov, Yu. A. Brychkov, O. I. Marichev**, *Integrals and Series,* Vol.1: Elementary Functions, Gordon and Breach 1986.

12.  **Andrei B. Shidlovskii**, *Transcendental Numbers,* W. de Grunter Studies in Math, 1989.

13.  **H. M. Srivastava and Junesang Choi**, *Series associated with the Zeta and related functions*, Kluwer, Boston, 2001.

14.  **Kenneth S. Williams, Zhang Nan-Yue**, Special Values of the Lerch Zeta Function and the Evaluation of Certain Integrals, *Proc. Amer. Math. Soc*., 119 (1), 35-49.

15.  **Eric W. Weisstein**. "Polylogarithm." From MathWorld. http://mathworld.wolfram.com/Polylogarithm.html